\documentstyle[12pt]{article}
\begin{document}
\newcommand{\n}{I{\hskip -5pt}N}
\newcommand{\r}{I{\hskip -5pt}R}
\newtheorem{definition}{Definition}[section]
\newtheorem{lemma}{Lemma}[section]
\newtheorem{theorem}{Theorem}[section]
\title{{Existence of Positive Solution of a Class of Semi-linear
Sub-elliptic Equation in the Entire Space $I{\hskip -5pt}H^n$
}
}
\author{{Zhujun Zheng}
\thanks{Email: zhengzj@mail.henu.edu.cn }\\
{\scriptsize Institute of Mathematics, Henan University, Kaifeng, 475001}\\
{Xiufang Feng}\\
{\scriptsize Department of Mathematics, Henan University, Kaifeng,
475001}
 }
\date{}
\maketitle
\begin{center}
\begin{minipage}{120mm}
\vskip 0.5cm {\bf Abstract}{~In this paper, we study the following
problem
$$
\left\{\begin{array}{ll}
\Delta_{H^n} u-u+u^p=0 ~~~~~~~~ & in H^n\\
u>0& in H^n\\
u(x)\rightarrow 0   &\rho(x)\rightarrow\infty
\end{array}
\right.
$$
where $1<p < \frac{Q+2}{Q-2}$, Q is the homogeneous dimension of
Heisenberg group $H^n$. Our main result is that this problem has
at least one positive solution.
}


{\bf AMS(1991) Subject Classification:~}35J60
\end{minipage}
\end{center}
\vskip 1cm
\baselineskip 20pt
\section{Introduction}
Let $H^n$ be the Heisenberg group, where $
\Delta_{H^n}=\sum\limits_{i=1}^n(X_i^2 + Y_i^2) $ is its
subelliptic Laplacian operator,
 $\rho (x)$ is the distance function from $x$ to the point $0$.
Under the real coordinate $(x_1,\cdots,x_n,y_1,\cdots,y_n,t)$, the
vector field $X_i$ and $Y_i$ are defined by
$$
\begin{array}{ll}
X_i=\frac{\partial}{\partial x_i}+2y_i\frac{\partial}{\partial t}&\\
&i=1,2, \cdots,n.\\
Y_i=\frac{\partial}{\partial y_i}-2x_i\frac{\partial}{\partial t}&\\
\end{array}
$$
and the distance function $\rho (x)$ is defined by
$$
\rho (x)=(\sum\limits_{i=1}^{n}(x_i^2+y_i^2)^2+t^2)^{\frac{1}{4}}.
$$

It is well known that $\{X_i, Y_i\}$ generate the real Lie algebra
of Lie group $H^n$ and
$$
[X_i, Y_i]=4\delta_{ij}\frac{\partial}{\partial t},~~
i,j=1,\cdots,n.
$$
In this Lie group, there is a group of natural  dilations defined
by
$$
\delta_\lambda(x,y,t)=(\lambda x, \lambda y, \lambda^2
t),~~\lambda >0
$$
where $x=(x_1,\cdots,x_n), y=(y_1,\cdots,y_n)$. With this group of
dilations, the Lie group $H^n$
 is a two step stratified nilpotent Lie group of homogeneous dimension $Q=2n+2$,
 and $\Delta_{H^n}$ is homogeneous partial differential operator of degree 2.
In this paper, we deal with the existence of the positive solution
to the following semi-linear subelliptic equation
\begin{equation}
\label{e1}
\left\{
\begin{array}{ll}
\Delta_{H^n} u-u+u^p=0 ~~~~~~~~ & in H^n\\
u>0& in H^n\\
u(x)\rightarrow 0   &\rho(x)\rightarrow\infty
\end{array}\right.
\end{equation}
where $1<p\leq \frac{Q+2}{Q-2}$

Equation (\ref{e1}) comes from the CR-Yamabe problem(see [14]) and
has been studied by several authors(see [4], [10],[12] and the
references therein). In the paper [12], they studied the problem
$$
\left\{
\begin{array}{ll}
\Delta_{H^n} u+u^p=0 ~~~~~~~~ & in H^n\\
u>0& in H^n\\
\end{array}\right.
$$
and showed that if the problem's solution is cylindrical, then it
must be 0. In the works [2] and [4], for $1<p<\frac{Q+2}{Q-2}$
they have gotten some results on the existence of the boundary
value problem of equation (\ref{e1}) on the bounded domain and
unbounded domain with thin condition. In these condition, the
corresponding functional satisfies P.S condition, and the normal
variational methods works. In the entire space $H^n,
1<p\leq\frac{Q+2}{Q-2}$, the  Folland-Stein-Soblev embedding lost
compactness. The corresponding functional lost P.S condition.
Their methods don't work. To our knowledge, in these situation,
there exists no report of progress on this problem up to now.

In the Euclidean space, the similar problem was studied by many
peoples(see [2], [3],[11], and the references therein). In [2],
W-Y Ding and W-M Ni gave some beautiful results on the similar
semilinear problem in Euclidean. But for our problem, as a
consequence of [12], it may not have radical symmetry solution. So
our problem is more subtle then them. This is one of the reasons
to study the problem (1).






 Our main result is the following theorem.

{\bf Theorem 1~} For $1<p< \frac{Q+2}{Q-2}$, the problem
(\ref{e1}) has a solution $u\in E$.


To proof our theorem, we first give some preliminary definition
and Lemmas. For $u\in C_0^\infty(H^n)$ the $C^\infty$ smooth
funciotn with compact support, we define norm $\| \|$ by
\begin{equation}
\label{e2} \|u\|^2=\int\limits_{\r^n}|\nabla_Hu|^2+u^2
\end{equation}
where
$\nabla_H=(\nabla_{X_1},\cdots,\nabla_{X_2},\nabla_{Y_1},\cdots,\nabla_{Y_n})$.
 Then we define the Folland-Stein-Sobolev space by
$E=\overline{C^\infty_0(H^n)}$, the completion of
$C_0^\infty(H^n)$ under the norm (\ref{e2}). This is a Hilbert
space. For
 $\Omega\subset H^n$,
the completion of $C_0^\infty(\Omega)$ in $E$ is denoted by
$E(\Omega)$, it is a Hilbert space too.  These spaces have
embedding theorem like the Sobolev embedding.

{\bf Lemma 1.1~} $\forall u\in E, 1<q\leq\frac{2Q}{Q-2}$, we have
\begin{equation}
\label{e3} \|u\|_{L^q}\leq C\|u\|
\end{equation}
where $C$ is a constant independent of $u$.

{\bf Lemma 1.2~} Let $\Omega\subset H^n$ be bounded smooth domain
in $H^n$, the embedding
\begin{equation}
\label{e4} E(\Omega)\hookrightarrow L^p(\Omega), ~~~1\leq
p<\frac{2Q}{Q-2}
\end{equation}
is compact.

We use two methods to solve the problem (\ref{e1}).

The first method:

In the Folland-Stein-Sobolev space $E$, we define the energy functional
\begin{equation}
\label{e5}
J(u)=\frac{1}{2}\int\limits_{H^n}|\nabla_Hu|^2+u^2-\frac{1}{p+1}\int u^{p+1}, ~~u\in E
\end{equation}
Let $B_k$ be the ball $B_k=\{x\in H^n\Big| \rho(x)<k\}$. Denoted
the completion $C_0^\infty(B_k)$ in $E$ by $E_k$, then
$$
\begin{array}{l}
E_k\subset E_{k+1}\subseteq E\\
\\
E=\overline{\bigcup\limits_{k=1}^\infty E_k}
\end{array}
$$
Set $J_k=J\Big| _{E_k}$, choose an element $u_0\in E_1\subseteq
E_1\subseteq\cdots\subseteq E_k\subseteq E$, such that
\begin{equation}
\label{e6}
J(u_0)<0, ~~~J_k(u_0)<0
\end{equation}
Let $\Gamma, \Gamma_k$ defined by
\begin{equation}
\label{e7}
\begin{array}{l}
\Gamma=\{r:[0,1]\rightarrow E\Big| r(0)=0,~ r(1)=u_0,~ r~ \hbox{is continuous}\}\\
\Gamma_k=\{r:[0,1]\rightarrow E_k\Big| r(0)=0,~ r(1)=u_0,~ r~ \hbox{is continuous}\}
\end{array}
\end{equation}
 Define

\begin{equation}
\label{e8}
\begin{array}{l}
c=\min\limits_{r\in \Gamma}\max\limits_{0\leq t\leq 1} I(r(t)),\\
c_k=\min\limits_{r\in \Gamma_k}\max\limits_{0\leq t\leq 1} I_k(r(t)).
\end{array}
\end{equation}
For $\Gamma_k\subseteq\Gamma_{k+1}\subset\Gamma$, we have
\begin{equation}
\label{e9}
c_k\geq c_{k+1}\geq c>0
\end{equation}

By mountain-path Lemma, we know $c_k$ is a critical value of the
functional $I_k$. Let $u_k$ be a critical point of $I_k$
corresponding the critical value, that is $I_k(u_k)=c_k$ and
$I'_k(u_k)=0$. By some complex estimates of $u_k$, we shall proof
$c$ is a critical value of $I$, and $u_k\rightarrow u$ in $E, u$
is a critical point and $I(u)=c$. By the maximum principle we get
a positive solution of (\ref{e1}).

the second method:

Define $M=\{u\in E\Big|\int\limits_{H^n}|u|^{p+1}=1\}\subset E$.
On the manifold, we define
\begin{equation}
\label{e10}
I(u)=\frac{1}{2}\int\limits_{H^n}|\nabla_H u|^2+u^2
\end{equation}

The main idea is, for $I$ have bounded from below, we define
\begin{equation}
\label{e11}
c=\inf\limits_{u\in M}I(u)
\end{equation}
Then we prove that the critical $c$ can be arrived by $u\in M$.
Then by Lagrane multiplier method, we know the problem have a
positive solution.



For the second method of subcritical case
to overcome the difficult that the functional $I$ lost P.S
condition, we give some Lions' version concentration-compactness
Lemmas. This is one of bones in this work.

For the case $1<p<\frac{Q+2}{Q-2}$, by our proof we know, for
every smooth bounded domain $\Omega$, the Dirichlet problem
\begin{equation}
\label{e13}
\left\{\begin{array}{ll}
\varepsilon^2\Delta_{H^n}u-u+u^p=0~~~~& in~ \Omega\\
u>0 & in ~\Omega\\
u=0 & on ~\partial\Omega
\end{array}
\right.
\end{equation}
have a least energy solution $u_\varepsilon$. Let
$x_\varepsilon\in\Omega, u(x_\varepsilon)=\max\limits_{x\in
\Omega} u(x)$, like J.Wei in the paper [13], we want to know what
is the $\lim\limits_{\varepsilon\rightarrow 0}{\rm
dist}(x_\varepsilon, \partial\Omega)$. In one of our preparing
works ([1]), we shall proof that
$$
{\rm dist}(x_\varepsilon, \partial\Omega)\rightarrow\max\limits_{x\in\Omega}d(x,\partial\Omega), ~~\varepsilon\rightarrow 0
$$
and we shall publish this result elsewhere. In the following, as
the Euclidean case, we shall study the effect of topology of the
unbounded $\Omega$.

\section{The Proofs of Main Theorem }

{\bf 2.1 The first method}

In this subsection, we shall use the mountain-path lemma and
domain extension method to proof the Theorem in the subcritical
exponent case $1<p<\frac{Q+2}{Q-2}$. And more, we get that the
problem have a least energy solution, and proof that
$$
c=\inf\limits_{r\in\Gamma}\max\limits_{0\leq t\leq 1} J(r t)
$$
can be arrived by a path $r_0\in\Gamma$. This is the foundation of our paper [1].

For the Folland-Stein-Sobolev embedding $E_k\hookrightarrow
L^{p+1}, 1<p<\frac{Q+2}{q-2}$  is compact ([5]), by the standard
method we have the following lemma.

{\bf Lemma 2.1~}
For $k\in\n$, the functional $J_k$ defined in the Hilbert $E_k$ satisfies P.S condition.

For an element $e\in E_1\subset E_k\subset E, \|e\|=1, \forall k\in\n$, we have
\begin{equation}
\label{e14}
J_k(te)=\frac{t^2}{2}-\frac{t^{p+1}}{p+1}\int|e|^{p+1}dx
\end{equation}
For $p+1>2$, we have the following Lemma 2.2.

{\bf Lemma~ 2.2~} For any $k\in\n $, there exists an element
 $u\in (\bigcap\limits_{k=1}^\infty E_k)\cap E$,
 such that
\begin{equation}
\label{e15}
I_k(u_0)<0, ~~J(u_0)<0
\end{equation}

For $\|u\|=1$, we have
\begin{equation}
\label{e16}
J(tu)=\frac{t^2}{2}-\frac{t^{p+1}}{p+1}\int\limits_{H^n}|u|^{p+1}
\end{equation}
By the Lemma 2.1, there is a positive constant $C>0$ independent of $u$, such that
\begin{equation}
\label{e17}
\int\limits_{H^n}|u|^{p+1}\leq C
\end{equation}

Combine  the inequality (\ref{e17}) and the formula (\ref{e16}) we have
\begin{equation}
\label{e18}
J(tu)\geq\frac{t^2}{2}-\frac{t^{p+1}}{p+1}C
\end{equation}
Since $p+1>2$, we have the following Lemma 2.3.

{\bf Lemma 2.3~} There is a neighborhood $U_k$ of 0 respectively
in $E_k$, and a neighborhood $U$ of 0 in $E$, such that
\begin{equation}
\label{e19}
J_k(u)\geq \alpha, ~~J(u)\geq\alpha
\end{equation}
for all $u\in U_k$ or $u\in U$ respectively, where $\alpha>0$ is a
positive constant.

From mountain path lemma and the above Lemma, we have the
following Lemma 2.4.

{\bf Lemma 2.4~} The value $c_k$ is a critical value of functional
$I_k$, and more we have
\begin{equation}
\label{e20}
c_k\geq c_{k+}\geq c>\alpha >0
\end{equation}

 Suppose $u_k$ is a critical point of $J_k$ corresponding the critical value $c_k$. Then we have
\begin{equation}
\label{e21}
J'(u_k)u_k=\|u_k\|^2-\int\limits_{H^n}|u|^{p+1}
\end{equation}

\begin{equation}
\label{e22}
J(u_k)=\|u_k\|^2-\frac{1}{p+1}\int\limits_{H^n}|u_k|^{p+1}=c_k>\alpha
\end{equation}

From (\ref{e21}) and (\ref{e22}), we have
\begin{equation}
\label{e23}
c_1\geq \frac{p}{p+1}\|u_k\|^2=c_k>\alpha
\end{equation}

That is to say $\{u_k\}$ is a bounded point set in $E$. So there
is a subset of $\{u_k\}$, we still denote it by $\{u_k\}$, and a
point $\overline u\in E$, such that

\begin{equation}
\label{eee}
u_k \rightharpoonup \overline u,
\end{equation}
and
$\overline u\geq 0$ is a weak solution of
$$
\Delta_{H^n} u-u+u^p=0 ~~~~~~~~  in H^n.
$$

By the method of Ding and Ni(see [2]), If we can prove $\overline
u\not\equiv 0$, then $\overline u$ is a critical of functional
$J$, and
$$
J(\overline u)=c.
$$
Then by maximum principal, we know $\overline u$ is a positive
solution of problem (\ref{e1}), and it is a positive least energy
solution of it. So if we can prove $\overline u \not\equiv 0$, our
theorem is proved. Next we locus on this problem.

For $u_k\in E$ is a solution of
$$
\Delta_{H^n} u-u+u^p=0 ~~~~~~~~  in H^n,
$$
we have
$$
\int|\nabla_{H^n}u_k|^2+u_k^2-\int u_k^{p+1}=0.
$$
Then we have
$$
\int u_k^2(u_k^{p-1}-1)=\int|\nabla_{H^n}u_k|^2\geq 0.
$$
Since $u_k\not\equiv 0$, there must be exists $\xi_k\in H^n$, such
that
\begin{equation}
\label{ee1}
u_k(\xi_k)=\max_{H^n}u_k\geq 1
\end{equation}
We claim that $\{\xi_k\}$ is a bounded subset of $H^n$. This is
our next lemma.

{\bf Lemma 2.5} The subset $\{\xi_k\}$ defined by (\ref{ee1}) is a
bounded subset of $H^n$.

Proof. For $\{u_k\}$ is bounded subset of $E$, by some standard
estimates and the Folland-Stein-Sobelev embedding theorem, there
is a positive constant $\alpha$, such that
$$
\sup_{H^n}u_k\leq \alpha.
$$
So there is a large enough $\beta > 0$ such that
\begin{equation}
\label{ee2}
-\Delta_H u_k+\beta u_k= u^p_k - (\beta - 1)u_k \leq 0
\end{equation}

Define function $v=ce^{-\delta \rho(x)},$ on $H^n$ where $c$ and
$\delta$ are positive number which shall be determined.

For $\Delta_H$ is a 2 order operator. So $\triangle_H v$ is a -1
order function. Then there are large positive numbers $R_0$, and
$\delta>0$, such that for all $\xi, \rho(\xi)> R_0$, and large
positive number $\beta '$ such that
\begin{equation}
\label{ee3}
-\Delta_H v(\xi)+ \beta ' v(\xi) \geq 0.
\end{equation}

Choose large positive number $R_0$, for all $\xi, \rho(\xi)=R_0,$
we have
\begin{equation}
\label{ee4}
(v-u)(\xi)\geq 0
\end{equation}
Set $\beta ''= \max\{\beta, \beta '\},$, then by (\ref{ee2},
\ref{ee3}, \ref{ee4})we have

$$
\left
\{\begin{array}{ll}
-\Delta_{H^n} (v-u)+\beta ''(v-u)\geq 0~~~~~~~~ & in H^n\backslash B_{R_0}(0),\\
v-u\geq 0, ~~~~~~~on \partial(H^n \backslash B_{R_0}(0)
\end{array}
\right.
$$

By the maximum principle, this implies that for all $\xi> R_0,$ ,
for any $k$,

\begin{equation}
\label{ee5}
u_k \leq ce^{-\delta\rho(\xi)}
\end{equation}

The inequality implies that ${\xi_k}$ is bounded.

For the Folland-Stein-Sobolve spaces have similar embedding
theorems with the Sobolev embedding and the Sub-Laplacian operator
have similar characters with the Laplacian operator(see [5]), so
by the method of Noussair, Ezzat S. and Swanson, Charles A(see
[13]), we have the following lemma.

{\bf Lemma 2.6.} There is a subsequence of $\{u_k\}$ we still
denote it by $\{u_k\}$, such that for any bounded domain $\Omega$,
$u_k \rightarrow \overline{u}$ in $C^{2+\alpha}(\Omega)$, where
$\alpha$ is a positive number. That is $u_k \rightarrow
\overline{u}$ in $C^{2+\alpha}_{loc}(H^n)$.

From Lemma 2.5 and Lemma 2.6, we have the following Lemma.

{\bf Lemma 2.7.} The functional defined by (\ref{eee})
$$
{\overline{u}\not\equiv 0}.
$$

Proof. For $\{\xi_k\}$ is bounded in $H^n$, so we may assume that
there is a $\xi_0\in H^k$, such that $\xi_k\rightarrow \xi_0$. So
we have
$$
u_k(\xi_k)\rightarrow \overline{u}(\xi_0).
$$
By the inequality (\ref{ee1}), we have $u(\xi_0)\geq 1$. That is
to say $u\not\equiv 0$. ~~~~~$\#$

{\bf 2.2 the second method}

In this subsection, we shall use the constraint functional method
to study the problem (1). First we defined the manifold
\begin{equation}
\label{e23.1}
M=\{u\in E\Big| \int|u|^{p+1}dx=1\}
\end{equation}
On this manifold, define a functional
\begin{equation}
\label{e24} I(u)=\frac{1}{2}\int|\nabla_Hu|^2+u^2, ~~\forall u\in
M
\end{equation}
It is obviously that the functional $I$ is bounded from below. We
shall study whether the functional defined by (\ref{e24}) arrive
its minimum on the manifold $M$. That is we want to find a $u_0\in
M$, such that
\begin{equation}
 \label{e25}
I(u_0)=\min\limits_{u\in M}I(u)=\alpha
\end{equation}

For the embedding $E\hookrightarrow L^{p+1}(H^n)$ lost
compactness, so the functional $I$ does not satisfy P.S condition.
To overcome this difficult, we first transplant  the Lions'
concentration-compactness Lemma([6,7,11]) to Heisenberg group
case.

{\bf Lemma 2.2.1~} Let $(\rho_m)_{m\geq 1}$ be a sequence in $L^1(H^n)$ satisfying:

\begin{equation}
\label{e25.1}
\rho_m\geq 0~ {\hbox{in}}~ H^n, ~\int\limits_{H^n}\rho_m=1
\end{equation}
Then there exists a sequence $(\rho_{m_k})_{k\geq 1}$ satisfying one the following three possibilities:

(i) (Compactness) ~There exists a sequence $z_k\in H^n$ such that
$\rho_{n_k}(z)$ is tight, i.e
\begin{equation}
\label{e26} \forall\varepsilon >0, \exists R<\infty,
\int\limits_{z_k+B_R}\rho_{n_k}(z)dz\geq 1-\varepsilon;
\end{equation}

(ii) (Vanishing)
$~\lim\limits_{k\rightarrow\infty}\sup\limits_{y\in
B_R}\int\rho_{n_k}(z)dz=0,$, for all $R<\infty$;

(iii) (Dichotomy) ~There exists $\alpha\in (0,1)$ such that for
all $\varepsilon>0$,  there $k_0\geq 1$ and $\rho_k^1, \rho_k^2\in
L^1_+(H^n)$ satisfying for $k\geq k_0$,
\begin{equation}
\label{e27}
\begin{array}{l}
\|\rho_{n_k}-(\rho_k^1+\rho^2_k)\|_{L^1}\leq \varepsilon\\
\\
|\int\limits_{H^n}\rho^1_kdz-\alpha|\leq\varepsilon
\end{array}
\end{equation}
and dist(supp$\rho_k^1, {\hbox{supp}}\rho_k^2)\rightarrow +\infty, k\rightarrow +\infty$, where $dz=dxdydt$.

For the measure $dxdydt$ on $H^n$, it has translation invariant
and it is a homogeneous on dilations $\delta_\lambda$ like the
measure on $\r^{2n+1}$. That is for $u\in L^1(H^n), z_0\in H^n$,
\begin{equation}
\label{e28}
\begin{array}{l}
\int\limits_{H^n}u(z)dz=\int\limits_{H^n}u(z\cdot z_0^{-1})dz\\
\int\limits_{H^n}u(z_\lambda z)dz=\lambda^{-Q}\int\limits_{H^n}u(z)dz
\end{array}
\end{equation}
Where $\lambda>0$. So, just like P.L.Lions[6,7], we can prove this
lemma. We omit its proof here.

Let $\{u_m\}\subset M, I(u)\rightarrow \min\limits_{u\in
M}I(u)=\alpha, m\rightarrow\infty$.  By the Folland-Stein-Sobolev
and there exists a constant $c>0$, such that
\begin{equation}
\label{e28.1}
\|u\|_p\leq c\|u\|, ~~\forall u\in E
\end{equation}
So $\alpha=\min\limits_{u\in M}I(u)>0$.

{\bf Lemma 2.2.2~} For the sequence $\{u_m\}$, there is a positive number $\{R_m\}$, for the function
\begin{equation}
\label{e29}
\nu_m(z)=R_m^{-\frac{1}{q}}u_m(\delta_{\frac{1}{R_m}}(z)
\end{equation}
 such that
\begin{equation}
\label{e30} \sup\limits_{z\in H^n}\int\limits_{B_1
(z)}|\nu_m|^q(w)dw=\frac{1}{2}=\int\limits_{B_1(0)}|\nu_m|^qdw
\end{equation}

Proof. For $u_m\in\{u_m\}, r>0, z_m^r\in H^n$, we define
\begin{equation}
\label{e31}
u_m^r=r^{n/q}u_m(\delta_{\frac{1}{r}}(z\cdot z_m^r))
\end{equation}
From (\ref{e28}), we have
\begin{equation}
\label{e33}
\int\limits_{H^n}\|u_m^r\|^q=r^{-n}\int\limits_{H^n}|u_m(\delta_{\frac{1}{r}}(z
z_m^r))|^q=\int\limits_{H^n}|u_m|^q=1
\end{equation}
So there exists a $R_m$, for every $z'_m\in H_n$,
\begin{equation}
\label{e34}
\int\limits_{B_1(z'_m)}|u^{p+1}_m|^qdz=\int\limits_{B_{R_m}(0)}|u_m|^qdz=\frac{1}{2}
\end{equation}
Define $\nu_m(z)=R_m^{-2/n}u_m(\delta_{\frac{1}{R_m}}z)$. From the formula (\ref{e34}), we have
$$
\sup\limits_{z\in H^n}\int\limits_{B_1(z)}|\nu_m|^qdx=\int\limits_{B_1(0)}|\nu_m|^qdx=\frac{1}{2}~~~~~~~~\#.
$$

Let $\rho_m=|\nu_m|^q$, then $\rho_m\in L^1(H^n)$, and
$\int\limits_{H^n}\rho_m=1$. From Lemma 2.2.2, we know case (ii)
in Lemma 2.2.1  can't occurs. We declare that the case (iii) can't
also. That is our following lemma.

{\bf Lemma 2.2.3}~ For the function $\rho_m\in L^1(H^n)$ defined
above, there is $z_m\in H^n$, such that $\rho_m(z\cdot z_m^{-1})$
is tight, i.e. there exists a number $R>0$ large enough, such that
\begin{equation}
\label{e35}
\int\limits_{z_m\cdot B_R(0)}\rho_m(z)dz\geq 1-\varepsilon
\end{equation}

Proof. By th Lemma 2.2.1 and Lemma 2.2.2, we only need prove  the
case(iii) in Lemma 2.2.1 does't occur. On contrary, there is a
number $\beta\in(0,\lambda)$ such that for all $\varepsilon>0$,
there exist $m_0\geq 1$ and $\rho^1_m, \rho_m^2\in L^1(H^n)$
satisfies for $m>m_0$,
\begin{equation}
\label{e36}
\begin{array}{l}
\|\rho_m-(\rho_m^1+\rho_m^2)\|_m \leq \varepsilon\\
\\
|\int\limits_{H^n}\rho_m^1 dz-\beta|\leq\varepsilon\\
\\
|\int\limits_{H^n}\int\rho^2_mdz-(1-\beta)|\leq\varepsilon
\end{array}
\end{equation}
and ${\hbox {dist(supp}}\rho_m^1, {\hbox {supp}}\rho_m^2)\rightarrow +\infty$.

Choose $r_m>0$, such that supp$\rho_m^1\subset B_{r_m}(0)$,
supp$\rho_m^2\subset H^n\backslash B_{r_m}(0)$, and
$r_m\rightarrow +\infty$ as $m\rightarrow +\infty$. Set
$\varphi\in C_0^\infty (B_2(0))$ such that $\varphi\equiv 1$ in
$B_1(0), 1\leq \varphi\leq 1$ and let $\varphi_m(\frac{x}{r_m})$.
Decompose
$$
\nu_m=\varphi_m \nu_m+(1-\varphi_m)\nu_m
$$
Then
\begin{equation}
\label{e37}
\begin{array}{ll}
\int\limits_{H^n}| \nabla_H \nu_m|^2+|\nu_m|^2=&\int\limits_{H^n}|
\nabla_H(\varphi_m\nu_m)|^2+\int\limits_{H^n}(\varphi_m\nu_m)^2+\int\limits_{H^n}|\nabla_H(1-\varphi_m)\nu_m)|^2\\
&~+\int\limits_{H^n}(1-\varphi_m)\nu_m)^2+2\int\limits_{H^n}\nabla_{H^n}(\varphi_m\nu_m)\cdot (1-\varphi_m)\nu_m)\\
&~+2\int\limits_{H^n}\varphi_m\nu_m(1-\varphi_m)\nu_m
\end{array}
\end{equation}
Next we estimate the last two terms in formula (\ref{e37}) respectively.
$$
\begin{array}{l}
\int\limits_{H^n}\nabla_H(\varphi_m\nu_m)\cdot\nabla_H((1-\varphi_m)\nu)\\
\geq -\int\limits_{H^n}|\nabla_H(\varphi_m\nu_m)\cdot\nabla_H(1-\varphi_m)\nu_m)|\\
\\
\geq -\int\limits_{H^n}|\nabla_H(\varphi_m\nu_m)| |\nabla_H((1-\varphi_m)\nu_m))|\\
\\
=-\int\limits_{B_{2r_m}(0)\backslash B_{r_m}(0)}|\nabla_H(\varphi_m\nu_m)| |\nabla_H((1-\varphi_m)\nu_m)|\\
\\
\geq\frac{1}{2}\Big[\int\limits_{B_{2r_m}(0)\backslash B_{r_m}(0)}|\nabla_H(\varphi_m\nu_m)|^2+\int\limits_{B_{2r_m}(0)\backslash B_{r_m}(0)}|\nabla_H((1-\varphi_m)\nu_m)|^2\Big]\\
\\
=-\frac{1}{2}[\int\limits_{B_{2r_m}(0)\backslash B_{r_m}(0)}\Big\{|\nabla_H(\varphi_m|^2\nu_m^2+2\nabla_H\varphi_m\cdot\nabla_H\nu_m\cdot\varphi_m\nu_m+\varphi_m^2|\nabla_m\nu_m|^2)|\\
~~+|\nabla_h\varphi_m|^2\nu_m^2-2\nabla_H\varphi_m\cdot\nabla_H\nu_m\cdot\varphi_m\nu_m+(1-\varphi_m)^2|\nabla_m\nu_m|^2\Big\}\\
\\
\geq -c\int\limits_{B_{2r_m}(0)\backslash B_{r_m}(0)}\nu_m^2+|\nabla_H\nu_m|^2\\
\\
\geq c\int\limits_{B_{2r_m}(0)\backslash B_{r_m}(0)}|\nu_m|^{p+1}
\end{array}
$$
Then we have
\begin{equation}
\label{e38}
\int\limits_{H^n}\nabla_H(\varphi_m\nu_m)\cdot((1-\varphi_m)\nu_m)+\int\varphi_m\nu_m(1-\varphi_m)\nu_m\geq -c\int\limits_{B_{2r_m}(0)\backslash B_{r_m}(0)}|\nu_m|^{p+1}
\end{equation}
From the proof of Lemma 2.2.1, we have $\forall \varepsilon>0, \exists m_0$, such that for $m>m_0$,
\begin{equation}
\label{e40}
\begin{array}{l}
\int\limits_{B_{r_m}(0)}|\nu_m|^{p+1}\leq \beta\\
\\
\int\limits_{H^n\backslash B_{r_m}(0)}|\nu_m|^{p+1}\leq 1-\beta+\varepsilon
\end{array}
\end{equation}
So by the inequalities of (\ref{e36}) and (\ref{e40}), we have
\begin{equation}
\label{e41} \int\limits_{B_{2r_m}(0)\backslash
B_{r_m}(0)}|\nu_m|^{p+1} \leq
c[\int\limits_{H^n}|\nu_m|^{p+1}-\int\limits_{H^n}(\rho^1_m+\rho+m^2)]+\varepsilon
\end{equation}
The inequality (\ref{e41}) means that
\begin{equation}
\label{e42}
\int\limits_{B_{2r_m}(0)\backslash B_{r_m}(0)}|\nu_m|^{p+1}=o(1)
\end{equation}
where $o(1)\rightarrow 0, m\rightarrow +\infty$.

Combine the formula (\ref{e42}),(\ref{e37}) and the inequality (\ref{e38}) we have
\begin{equation}
\label{e43}
\int\limits_{H^n}|\nabla_H\nu_m|^2+|\nu_m|^2=\|\varphi_m\nu_m\|^2+\|(1-\varphi_m)\nu_m\|^2+o(1)
\end{equation}

By the Folland-Stein-Sobolev embedding and formula (\ref{e43}), we have
\begin{equation}
\label{e44}
\begin{array}{ll}
\|\nu_m\|^2&=\|\varphi_m\nu_m\|^2+\|(1-\varphi_m)\nu_m\|^2+o(1)\\
\\
&\geq S(\|\varphi_m\nu_m\|_{L^{p+1}}^{\frac{2}{p+1}}+\|(1-\varphi_m)\nu_m\|^{\frac{2}{p+1}})+o(1)\\
\\
&\geq S(\int\limits_{H^n}\rho_m^1)^{\frac{1}{p+1}}+(\int\limits_{H^n}\rho_m^2)^{\frac{2}{p+2}})+o(1)\\
\\
&\geq S(\beta^{\frac{2}{p+1}}+(1-\beta)^{\frac{2}{p+1}})+o(1)
\end{array}
\end{equation}

By the define of $\alpha $ and the independence of domain of the
best Folland-Stein-Sobolev constant we know $S=\alpha$. And by the
define $\nu_m$ we have
$\|\nu_m\|^2\rightarrow\alpha(m\rightarrow\infty)$. So by the
inequality (50)we have
\begin{equation}
\label{e45}
\alpha\geq\alpha(\beta^{\frac{2}{p+1}}+(1-\beta)^{\frac{2}{p+1}})
\end{equation}

For $\frac{2}{p+1}<1, 0<\beta<1$, we get $\alpha>\alpha$, that is
a contradiction. So the case (iii) of Lemma 2.2.1 can't occur.

{\bf Theorem 2.2.1}~ There is a subsequence of $\{\nu_m\}$, we
still  denote it by $\{\nu_m\}$, there exists a point $u_0\in M$,
such that $\nu_m\rightarrow\nu_0$ in $E$, and
$$
I(\nu_0)=\alpha
$$

Proof: For $\{\nu_m\}$ is bounded in $E$, we have  a subsequence
of it, and we still denote it by $\{\nu_m\}$, and there exists a
$\varepsilon, \nu_0\in E$ such that $\nu_m\rightarrow \nu_0$. For
$\varepsilon<\frac{1}{2}$, and Lemma 2.1.2, we get $(z_m\cdot
B_R(0))\cap B_1(0)\neq \emptyset$, so the points sequence
$\{z_m\}$ is a bounded set. That is implies that tere is a
subsequence of $\{z_m\}$, we still denote it by $\{z_m\}$ and a
point $z_0\in H^n$, such that $z_m\rightarrow z_0,
(m\rightarrow\infty)$. Then we have
\begin{equation}
\label{e46}
\int\limits_{z_0\cdot B_{1+2R}(0)}|\nu_m|^{p+1}>1-\varepsilon
\end{equation}

From the Folland-Stein-Sobolev emedding, we know there is a
subsequence $\{\nu_m\}$ we still denote it by $\{\nu_m\}$, such
that $\nu_m\rightarrow\nu_0, m\rightarrow \infty$ in
$H^{1,2}(B_{1+2R}(0))$ and
\begin{equation}
\label{e47}
\int\limits_{z_0\cdot B_{1+2R}(0)}|\nu_0|^{p+1}>1-\varepsilon
\end{equation}
From the Fatou Lemma we know
\begin{equation}
\label{e48} \int\limits_{H^n}|\nu_0|^{p+1}\leq
\lim\limits_{m\rightarrow\infty}\int\limits_{H^n}|\nu_m|^{p+1}=1
\end{equation}
Combine the inequalities of (\ref{e47}) and (\ref{e48}), we have
\begin{equation}
\label{e49}
\int\limits_{H^n}|\nu_0|^{p+1}=1
\end{equation}
This implies that $\nu_m\rightarrow \nu_0, m\rightarrow\infty $ in $E$. So we have
$$
I(\nu_0)=\lim\limits_{m\rightarrow\infty}I(\nu_m)=\alpha
$$

By the Lagrange multiplier, there is a positive number, such that
$$
\Delta_{H^n}\nu_0-\nu_0+\lambda \nu_0^p=0
$$
Set $u_0=\lambda^{p-1}u$, then we have
$$
\Delta_{H^n}u_0-u_0+u_0^{p}=0
$$
By the maximum principle, we get our our main theorem 1.1.

\vskip 1cm

\begin{center}
 {\large\bf Acknowledgments}
\end{center}
I would like thanks Prof. Zhouping Xin for his inviting me to
visit IMS and useful discussions. I would like to thanks Prof.
Juncheng Wei for he let me notice this problem and helpful
discussions. I would like thanks Prof. Changfeng Gui, Prof.
Yongsheng Li, Prof. Quanshen Jiu and Prof. Jiabao Su for their
useful discussions. This work is finished during my visiting IMS
of The Chinese University of Hong Kong. And is partially supported
by the Zheng Ge Ru Foundation, Mathematical Tianyuan Fund with
grand No. 10226002 and the Natural Science Foundation of
Educational Committee of Henan Province with grand Number
2000110010










\end{document}